\documentclass[11pt,leqno]{article}
\usepackage{amsfonts}
\newtheorem{theo}{Theorem}[section]
\newtheorem{lem}[theo]{Lemma}
\newtheorem{pro}[theo]{Proposition}
\newtheorem{co}[theo]{Corollary}

\def\Box{{\vcenter{\vbox{\hrule height.4pt
    \hbox{\vrule width.4pt height7pt \kern7pt   \vrule width.4pt}
        \hrule height.4pt}}}}

\title{\Large \bf {Vanishing theorems on Hermitian manifolds }}
\author{{\sc B.Alexandrov} \hspace{3mm} {\sc S.Ivanov}
\thanks{The authors are supported by Contract MM 809/1998 with the
Ministry of Science and Education of Bulgaria and by Contract 238/1998 with the
University of Sofia "St. Kl. Ohridski".}}
\date{}
\begin{document}
\maketitle \thispagestyle{empty} \vspace{2mm}
%\begin{center}
%Diff. Geom. Appl. {\bf 14} (3) (2001), 251-265.
% http://xxx.lanl.gov, math.DG/9901090,
%{\sc Department of Mathematics\\ University of Sofia\\
%"St. Kl. Ohridski" }
%\end{center}
\vspace{25mm}

\begin{abstract}
We prove the vanishing of the Dolbeault cohomology groups on Hermitian manifolds
with $dd^c$-harmonic K\"ahler form and positive $(1,1)$-part of the Ricci form
of the Bismut connection. This implies the vanishing of the Dolbeault cohomology
groups on complex surfaces which admit a conformal class of Hermitian metrics,
such that the Ricci tensor of the canonical Weyl structure is positive. As
a corollary we obtain that any such surface must be rational with $c_1^2>0$.
As an application, the pth
Dolbeault cohomology groups of a left-invariant complex structure compatible
with a bi-invariant metric on a compact even
dimensional Lie group are computed.
\\[15mm]
{\bf Running title:} Vanishing theorems on Hermitian manifolds
\\[5mm]
${\bf MS}$ {\bf classification: } 53C55; 53C15
\end{abstract}
\newpage

\section{Introduction and statement of the results}

\indent In \cite{H} Hitchin has proved that a
$Spin^c$-structure with determinant bundle $L$ on a K\"ahler manifold $M$ does not
admit harmonic spinors if the following positivity condition is satisfied: the
scalar curvature $s \ge 2\sum |\lambda _i|$, with strict equality at some
point. Here $\lambda _i$ are the eigenvalues
of the curvature form of the chosen unitary connection on $L$, considered as a
skew-symmetric endomorphism of the tangent bundle. As shown in \cite{H}, when
the $Spin^c$-structure is the canonical one and the Ricci tensor is positive,
the above result implies the vanishing of the Dolbeault cohomology
groups $H^p (M,{\mathcal O})$, thus rediscovering the Bochner-Kodaira vanishing
theorem.

On the other hand, any Hermitian manifold $(M,g,J)$ carries a unique Hermitian
connection with totally skew-symmetric torsion, the Bismut connection (cf.
\cite{B,G2}). A Hermitian manifold equipped with the Bismut
connection is also called K\"ahler with torsion or KT manifold, see e.g.
\cite{KT0}. KT manifolds arise in a natural way in physics as target
spaces of (2,0)-supersymmetric sigma models with Wess-Zumino term (torsion)
\cite{KT1,KT2} (see also \cite{homo} and the references there). A
characterization of KT manifolds in terms of the properties
of a two form is given in \cite{KT0}. KT structures on compact Lie groups and on
compact homogeneous spaces
are presented in \cite{SST} and \cite{homo} respectively.

The purpose of this note is to find conditions expressed in terms of the
Bismut connection, which imply the vanishing of the Dolbeault cohomology
groups on Hermitian manifolds. In Lemma~\ref{lem1} we give a slight
modification of the Lichnerowicz type formula for the Dolbeault operator,
proved by Bismut \cite{B}. As an application we obtain the following theorem:

\begin{theo}\label{th2}
Let $(M,g,J)$ be a compact $2n$-dimensional ($n>1$) Hermitian manifold with
K\"ahler form $\Omega $. Suppose that $\Omega $ is $dd^c$-harmonic, i.e.
$dd^c \Omega =0$ and $(dd^c)^*\Omega =0$. Suppose also that the $(1,1)$-part of
the Ricci form of the Bismut connection is non-negative everywhere on
$M$.

a) Then every $\overline{\partial}$-harmonic (0,p)-form, $p=1,2,...,n$, must be
parallel with respect to the Bismut connection.

b) If moreover the $(1,1)$-part of
the Ricci form of the Bismut connection
is  strictly positive at some point, then the cohomology groups
$H^p (M,{\mathcal O})$ vanish for $p=1,2,...,n$.
\end{theo}

The condition $(dd^c)^*\Omega =0$ is equivalent to $dd^c \Omega^{n-1} =0$ and
to $d^* \theta =0$, where
$\theta $ is the Lee form. As shown by Gauduchon \cite{G3}, any conformal
class of Hermitian
metrics on a compact manifold contains a unique (up to homothety) metric satisfying
$d^* \theta =0$. This metric is called the Gauduchon metric.

The condition $dd^c \Omega =0$ is well known.
The Hermitian manifolds with $dd^c \Omega =0$ are called strong KT manifolds
and arise as target
spaces of classical (2,0)-supersymmetric sigma models with
torsion  (see e.g.\cite{KT0}). Under the same condition Bismut \cite{B} has
proved a local index theorem
for the Dolbeault operator.

Any compact Lie group with bi-invariant metric and
compatible left-invariant complex structure has $dd^c$-harmonic K\"ahler form.
Moreover, such a group does not admit any K\"ahler
metric unless it is a torus. In section~\ref{sec5.1} we show how Theorem~\ref{th2}
can be applied to compute the Hodge numbers $h^{0,p}$ of a Lie group
with such a Hermitian structure.

Recall that any Hermitian manifold $(M,g,J)$ carries a
canonical Weyl structure, i.e. a torsion-free connection which preserves the
conformal class $[g]$ and depends only on $[g]$ and $J$ (cf. \cite{V2}). When
the dimension is 4 the canonical Weyl structure preserves the complex structure
$J$, but in higher dimensions this is true only for locally conformally
K\"ahler manifolds (cf. \cite{V2}).  The symmetric part of the Ricci tensor of
the canonical Weyl structure on Hermitian surface is of type $(1,1)$ (cf.
\cite{GI}). We notice that the corresponding $(1,1)$-form coincides with the
$(1,1)$-part
of the Ricci form of the Bismut connection (see Lemma~\ref{pro1} below).
Since in dimension 4 the conditions
$dd^c \Omega =0$ and $(dd^c)^*\Omega =0$ coincide, the existence of the
Gauduchon metric allows us
to restate Theorem~\ref{th2} in conformally invariant fashion:

\begin{theo}\label{th4}
Let $(M,J)$ be a compact complex surface. Suppose that there exists a conformal
structure $c$ compatible with $J$ and such that
the symmetric part of the Ricci
tensor of the canonical
Weyl structure is non-negative everywhere on $M$ and strictly positive at
some
point. Then the cohomology groups $H^p (M,{\mathcal O})$ vanish for $p=1,2$ and
$(M,J)$ must be a rational surface with $c_1 ^2 >0$.
\end{theo}

It is well known that ${\bf CP^1} \times {\bf CP^1}$, ${\bf CP^2}$ and the blow-ups of ${\bf CP^2}$ in up to 8 points in
general position  carry K\"ahler
metrics with positive Ricci tensor. On the other hand, the Einstein-Weyl
structures with positive Ricci tensor are in fact K\"ahler-Einstein, cf.
\cite{GI}.  In section~\ref{sec5.2} we show that any compact 4-dimensional K\"ahler manifold with positive Ricci form carries a non-K\"ahler Hermitian structure with $dd^c$-harmonic K\"ahler form and positive $(1,1)$-part of the Ricci form of the Bismut connection. These are also examples of manifolds satisfiying the assumptions of Theorem~\ref{th4} with a conformal class of Hermitian metrics not containing a K\"ahler metric. The product of such a 4-dimensional manifold with a compact $(2n-4)$-dimensional K\"ahler manifold with positive Ricci form is a compact $2n$-dimensional non-K\"ahler Hermitian manifold with $dd^c$-harmonic K\"ahler form and positive $(1,1)$-part of the Ricci form of the Bismut connection. Thus there are non-K\"ahler manifolds satisfying the assumptions of Theorem~\ref{th2} b) in any dimension. We are not aware of any manifolds which satisfy the assumptions of Theorem~\ref{th2} b) but do not admit any K\"ahler metric. It follows from Theorem~\ref{th4} that if such manifolds do exist, then their dimension must be greater than 4.

In the 4-dimensional case, i.e. when $n=2$, the
 condition $d\Omega = \theta \land \Omega$ is always satisfied.
 When $n>2$ this condition means that the Lee form $\theta$ is closed
and the manifold
is locally conformally K\"ahler. Note that for $n>2$ it is not possible to have
simultaneously $dd^c \Omega =0$ and $d\Omega = \theta \land \Omega$ on a
compact manifold, except in the K\"ahler case (see Remark~1
below). A particular subclass
of locally conformally K\"ahler manifolds is formed by the generalized Hopf manifolds, the non-K\"ahler Hermitian
manifolds whose Lee form is parallel with respect to the Levi-Civita connection. For locally conformally K\"ahler manifolds we prove

\begin{theo}\label{th3}
Let $(M,g,J)$ be a compact $2n$-dimensional ($n>1$) locally conformally
K\"ahler manifold and $g$ be the Gauduchon metric. Suppose that the $(1,1)$-part of the Ricci form of the Bismut
connection is non-negative on $M$. Then:

a) If $n$ is even and $n>2$, then the cohomology group
$H^{\frac{n}{2}} (M,{\mathcal O})$ vanishes or $(M,g,J)$ is a K\"ahler
manifold.

b) If (M,g,J) is a generalized Hopf manifold, then the first Betti
number $b_1=1$ and  for the Hodge numbers
we have
$$
h^{p,0}=0, \quad p=1,2,\ldots,n, \quad h^{0,q}=0, \quad q=2,3,\ldots,n,
\quad h^{0,1}=1.
$$
\end{theo}

{\it Acknowledgement:} We would like to thank V. Apostolov for his interest and
valuable suggestions on the present paper. We also thank to T.Pantev for  his
comments on the cohomology of compact Lie groups.

\section{Preliminaries}

Let $(M,g,J)$ be a $2n$-dimensional ($n>1$) Hermitian manifold with complex
structure $J$ and compatible metric $g$. Let $Spin^c(M)$ be the canonical
$Spin^c$-structure on $(M,g,J)$, i.e. the
$Spin^c$-structure whose determinant bundle is the anticanonical bundle
$K^{-1}$. Let $\Sigma $ be the space of spinors. Then we obtain the Clifford
module $Spin^c(M) \times _{Spin^c(2n)} \Sigma$ and as is well known it is
isomorphic to $\Lambda ^{0,\bullet }M$, the bundle of $(0,p)$-forms. Locally
we can chose a square root $K^{\frac{1}{2}}$ of $K$ and let $\Sigma M$ be the
bundle of spinors corresponding to $K^{\frac{1}{2}}$ (cf. \cite{H}). Thus
locally
\begin{equation}\label{1}
\Lambda ^{0,\bullet }M = \Sigma M \otimes K^{-\frac{1}{2}}.
\end{equation}
We shall denote the Clifford product of a form $\alpha \in  \Lambda ^\bullet M$
and $\psi \in \Lambda ^{0,\bullet }M$ by $\alpha \psi $.

Let $\Omega $ be the K\"ahler form of $(M,g,J)$, defined by
$\Omega (X,Y)=g(X,JY)$. Denote by $\theta$ the Lee form of $(M,g,J)$,
$\theta = \frac{1}{n-1} d^* \Omega \circ J$. For a 1-form $\alpha $ we
shall denote by $J\alpha $ the form dual to $J\alpha ^{\#}$, where
$\alpha ^{\#}$ is the vector dual to $\alpha $. Equivalently,
$J\alpha = -\alpha \circ J$. Hence, $d^* \Omega = (n-1)J\theta $.

The connections we shall use in the sequel are the Levi-Civita connection
$\nabla ^L$, the Chern connection $\nabla ^C$, the Bismut connection
$\nabla ^B$ and the Weyl connection $\nabla ^W$. Below we briefly recall some
of their properties, which we shall need.

The Chern connection is given by
\begin{equation}\label{2}
g(\nabla ^C _X Y,Z) = g(\nabla ^L _X Y,Z) + \frac{1}{2} d\Omega (JX,Y,Z).
\end{equation}
Restricted to $T^{1,0} M$ it coincides with the canonical connection of this
holomorphic bundle.

The Bismut connection is given by
\begin{equation}\label{3}
g(\nabla ^B _X Y,Z) = g(\nabla ^L _X Y,Z) + \frac{1}{2} d^c \Omega (X,Y,Z).
\end{equation}
Recall that $d^c = i(\overline {\partial} - \partial )$. In particular,
$d^c \Omega (X,Y,Z) = - d \Omega (JX,JY,JZ)$. This connection has been used by
Bismut in \cite{B} to prove a Lichnerowicz type formula for the Dolbeault
operator.
It is one of the canonical Hermitian connections (cf. \cite{G2}) and in the
set of all Hermitian connections it is characterized by the fact that it is the
only connection with totally skew-symmetric torsion.

The canonical Weyl connection determined by the Hermitian structure of $M$ is
the unique
torsion-free connection $\nabla ^W$ such that $\nabla ^W g = \theta \otimes g$.

The canonical Weyl connection is invariant under conformal changes of the
metric, since if $\widetilde{g} = e^f g$, then
$\widetilde{\theta }= \theta + df$. We have
\begin{equation}\label{4}
\nabla ^W_X Y = \nabla ^L_X Y - \frac{1}{2} \theta (X)Y
- \frac{1}{2} \theta (Y)X + \frac{1}{2} g(X,Y)\theta ^{\#}.
\end{equation}
We shall denote by $Ric^W (X,Y)$ the symmetric part of the Ricci tensor
$tr \{ Z \longrightarrow R^W (Z,X)Y \}$ of $\nabla ^W$ (it is easy to see that its skew-symmetric part is equal to $\frac{n}{2} d\theta $).

The canonical Weyl connection preserves the complex structure iff
\begin{equation}\label{101}
d\Omega = \theta \land \Omega
\end{equation}
(cf. \cite{V2}). The condition (\ref{101}) is always satisfied in dimension 4. In higher dimensions it means that the manifold is locally conformally K\"ahler and in particular, $d\theta =0$. A 4-dimensional manifold is locally conformally K\"ahler when $d\theta =0$.

The Chern and Bismut connections preserve the
Hermitian structure. Hence, they induce unitary
connections on $K^{-1}$ with curvatures $i\rho ^C$ and $i\rho ^B$, where
$$\rho ^C (X,Y)=\frac{1}{2}\sum_{j=1}^{2n} g(R^C (X,Y)e_j,Je_j)$$ is the
Ricci form of $\nabla ^C$, and the
Ricci form $\rho ^B$ of $\nabla ^B$ is
defined similarly. Here and henceforth $\{e_1,e_2,...,e_{2n}\}$ is a local
orthonormal frame
of the tangent bundle $TM$ and for the curvature we adopt the
following convention: $R(X,Y) = [\nabla _X,\nabla _Y] - \nabla _{[X,Y]}$.

By (\ref{2}) and (\ref{3}) we obtain
\begin{equation}\label{5}
\nabla ^C _X \varphi = \nabla ^B _X \varphi + (n-1)iJ\theta (X) \varphi, \qquad
\varphi \in \Gamma (K^{-1}).
\end{equation}
Hence,
\begin{equation}\label{6}
\rho ^C = \rho ^B + (n-1)dJ\theta .
\end{equation}

In the following we shall denote by $<.,.>$ and $|.|$ the pointwise inner
products and norms and by $(.,.)$ and $\| .\|$ the global ones
respectively.

The trace of $\rho ^C$, denoted by $2u$ in \cite{G4}, is defined by
$$2u=2<\rho ^C , \Omega >=\sum_{j=1}^{2n} \rho ^C (Je_j,e_j).$$ As proved in
\cite{G4},
\begin{equation}\label{7}
2u = s -(n-1)d^* \theta + \frac{1}{2} |d\Omega |^2,
\end{equation}
where $s$ is the scalar curvature of $g$. Let the trace of $\rho ^B$ be
$$b=2<\rho ^ B ,\Omega >=\sum_{j=1}^{2n} \rho ^B (Je_j,e_j).$$ By (\ref{6}) we
have
\begin{equation}\label{8}
2u=b + (n-1)\sum_{j=1}^{2n} dJ\theta (Je_j,e_j).
\end{equation}
But a direct computation yields
\begin{equation}\label{8.5}
\sum_{j=1}^{2n} dJ\theta (Je_j,e_j) = \frac{2}{n-1}<dd^* \Omega , \Omega > =
2(n-1) |\theta |^2 + 2d^* \theta .
\end{equation}
Hence, by (\ref{8}) and (\ref{8.5}) we obtain
\begin{equation}\label{110}
2u = b + 2(n-1)d^* \theta + 2(n-1)^2 |\theta |^2
\end{equation}
and it follows from this equality and (\ref{7}) that
\begin{equation}\label{9}
b = s - 3(n-1)d^* \theta - 2(n-1)^2 |\theta |^2+ \frac{1}{2} |d\Omega |^2 .
\end{equation}

In the sequel we shall need also the following equalities, which are
obtained by direct computations:
\begin{equation}\label{10}
<dd^c \Omega , \Omega \land \Omega > = 2(n-1)^2 |\theta |^2 - 2|d\Omega |^2
+ 2(n-1)d^* \theta ,
\end{equation}
\begin{equation}\label{11}
|\theta \land \Omega |^2 = (n-1) |\theta |^2 .
\end{equation}

\noindent {\bf Remark~1:} Recall that a Hermitian manifold is called balanced \cite{Mic} (or semi-K\"ahler, cf. e.g. \cite{Ga0}) iff its Lee form vanishes identically. It follows from (\ref{10}) that any balanced Hermitian manifold with $dd^c \Omega =0$ is K\"ahler. By (\ref{10}) and (\ref{11}) it is also clear that it is not possible to have
$dd^c \Omega =0$ and $d\Omega = \theta \land \Omega$ simultaneously on
compact Hermitian non-K\"ahler manifolds, except when $n=2$.

\vspace{3mm}
By choosing a metric connection on the tangent bundle $TM$ and a unitary
connection on the determinant bundle $K^{-1}$ we obtain a connection on
$Spin^c(M)$ and hence a unitary connection on $\Lambda ^{0,\bullet }M$. We
shall denote the connections obtained in this way by $\nabla $ with two upper
indices: the first denoting the connection on $TM$ and the second one  the
connection on $K^{-1}$. For example, $\nabla ^{L,C}$ is the connection obtained
from the Levi-Civita connection on $TM$ and the Chern connection on $K^{-1}$.

From (\ref{5}) and (\ref{3}) we obtain respectively
\begin{equation}\label{12}
\nabla ^C _X \varphi = \nabla ^B _X \varphi +
\frac{n-1}{2}iJ\theta (X) \varphi,
\qquad \varphi \in \Gamma (K^{-\frac{1}{2}}),
\end{equation}
\begin{equation}\label{13}
\nabla ^B _X \psi = \nabla ^L _X \psi +
\frac{1}{4} (\iota _X d^c \Omega )\psi , \qquad \psi \in \Gamma (\Sigma M).
\end{equation}
Here $\iota _X$ denotes the interior multiplication operator.

Using the local formulae (\ref{12}), (\ref{13}) and (\ref{1}), we get
\begin{equation}\label{14}
\nabla ^{B,C} _X \psi = \nabla ^{B,B} _X \psi + \frac{n-1}{2}iJ\theta (X) \psi,
\qquad  \psi \in \Gamma (\Lambda ^{0,\bullet }M),
\end{equation}
\begin{equation}\label{15}
\nabla ^{B,C} _X \psi = \nabla ^{L,C} _X \psi +
\frac{1}{4} (\iota _X d^c \Omega )\psi , \qquad
\psi \in \Gamma (\Lambda ^{0,\bullet }M).
\end{equation}

\section{Vanishing of the plurigenera}

Recall that for $m>0$ the $m$-th plurigenus of a compact complex manifold $M$ is defined by $p_m = \dim H^0 (M,{\cal O}(K^m))$.

The following result is a consequence of Gauduchon's plurigenera theorem \cite{Ga0,G4}.

\begin{pro}\label{tpr1}
Let $(M,g,J)$ be a compact $2n$-dimensional Hermitian manifold and $b$ be
the trace of
the Ricci form of the Bismut connection. If
\begin{equation}\label{201}
b \ge 0 ,
\end{equation}
then the plurigenera $p_m \le 1$ for all $m>0$. If furthermore the inequality (\ref{201}) is strict at some point or the Gauduchon metric of the Hermitian structure is not balanced, then $p_m =0$ for all $m>0$. The same conclusions hold when $g$ is the Gauduchon metric and in (\ref{201}) $b$ is replaced by $\int_M b\,dV$.
\end{pro}

\noindent {\it Proof:} Let $g$ be the Gauduchon metric. It follows from (\ref{110}) that if
\begin{equation}\label{202}
\int_M b\,dV \ge 0,
\end{equation}
then
\begin{equation}\label{203}
\int_M u\,dV \ge 0
\end{equation}
and the inequality (\ref{203}) is strict if (\ref{202}) is strict or $\theta $ is not identically zero. Hence, in this case the assertions of the proposition follow from Gauduchon's plurigenera theorem \cite{Ga0,G4}.
\\
It is well known (see e.g. \cite{Bes}) that under a conformal change of the metric $\widetilde{g} = e^f g$ the scalar curvarures of $\widetilde{g}$ and $g$ are related by
$$e^f \widetilde{s} = s - \frac{(2n-1)(2n-2)}{4}|df|^2 +(2n-1)\Delta f ,$$
where the norm and the Laplace operator are with respect to $g$. Using this equality and the fact that $d^*\theta =0$, we get the following relation between the traces of the Ricci forms of the Bismut connections of $\widetilde{g}$ and $g$:
$$e^f \widetilde{b} = b - (n-2)\Delta f - (n-1)(n-2)<df,\theta >.$$
Thus, if $\widetilde{b} \ge 0$, then we obtain (\ref{202}) and this proves the proposition. \hfill {\bf Q.E.D.}

\vspace{3mm}
The fact that if $b=0$, then $p_m \le 1$, is proved by Grantcharov \cite{Gr} by similar methods. He has applied it to HKT manifolds. Recall \cite{KT0} that a hyper-Hermitian manifold $M$ is called a HKT manifold if the Bismut connections of the three complex structures coincide (and hence the Bismut connections of the whole $S^2$-family of complex structures coincide). Hence, the common Bismut connection preserves the hyper-Hermitian structure and thus its holonomy is contained in $Sp(\frac{\dim M}{4})$. This implies in particular that its Ricci form is zero. Thus, applying Proposition~\ref{tpr1} we obtain that on a compact HKT manifold all the plurigenera of any complex structure of the $S^2$-family are less or equal to one.

\section{Vanishing of Hodge numbers}

In the following we shall denote by $\Box$ the Dolbeault operator
$\sqrt 2 (\overline {\partial} + {\overline {\partial}}^*)$ on
$\Lambda ^{0,\bullet }M$.

\begin{lem}\label{lem1}
Let $(M,g,J)$ be a compact $2n$-dimensional ($n>1$) Hermitian manifold. Then
for $\psi \in \Gamma (\Lambda ^{0,\bullet }M)$
\begin{eqnarray}\label{16}
\| \Box \psi \| ^2 = & & \| \nabla ^B \psi \| ^2 +
(n-1)\Re (iJ\theta \Box  \psi ,\psi ) \\& & + \frac{1}{4}
((b + 3(n-1)d^* \theta + (n-1)^2 |\theta |^2 - |d\Omega |^2) \psi ,\psi)
\nonumber  \\& & + \frac{i}{2}(\rho ^B \psi ,\psi ) +
\frac{1}{4}(dd^c \Omega \psi ,\psi). \nonumber
\end{eqnarray}
\end{lem}

\noindent {\it Proof:} Theorem~2.3 in \cite{B}, formulated in our notations,
yields
\begin{equation}\label{17}
\Box ^2 = (\nabla ^{B,C})^* \nabla ^{B,C} + \frac{s}{4} +
\frac{i}{2}\rho ^C + \frac{1}{4}dd^c \Omega - \frac{1}{8}|d^c \Omega |^2.
\end{equation}
Since $d^* J\theta = 0$, we have
$(\nabla ^{B,C} _{J\theta ^{\#}} \psi ,\psi ) +
(\psi ,\nabla ^{B,C} _{J\theta ^{\#}} \psi ) = 0$.
Hence, by (\ref{14}) we obtain
\begin{equation}\label{18}
\| \nabla ^{B,B} \psi \| ^2 = \| \nabla ^{B,C} \psi \| ^2 +
(n-1)i(\nabla ^{B,C} _{J\theta ^{\#}} \psi ,\psi ) +
\frac{(n-1)^2}{4}(|\theta |^2 \psi ,\psi).
\end{equation}
It follows from Theorem~2.2 in \cite{B} that
$$\Box = D^{L,C} + \frac{1}{4}d^c \Omega,$$
where $D^{L,C}$ is the Dirac operator of $\nabla ^{L,C}$,
$D^{L,C}\psi = \sum_{j=1}^{2n} e^j \nabla ^{L,C}_{e_j}\psi $. Hence,
\begin{equation}\label{19}
J\theta \circ \Box + \Box \circ J\theta =
J\theta \circ D^{L,C} + D^{L,C} \circ J\theta +
\frac{1}{4}J\theta \circ d^c \Omega + \frac{1}{4}d^c \Omega \circ J\theta .
\end{equation}
On the other hand,  we have
\begin{equation}\label{20}
J\theta \circ D^{L,C} + D^{L,C} \circ J\theta = dJ\theta + d^* J\theta -
2\nabla ^{L,C} _{J\theta ^{\#}} = dJ\theta - 2\nabla ^{L,C} _{J\theta ^{\#}},
\end{equation}
\begin{equation}\label{21}
J\theta \circ d^c \Omega + d^c \Omega \circ J\theta =
-2 \iota _{J\theta ^{\#}} d^c \Omega.
\end{equation}
Substituting (\ref{20}) and (\ref{21}) in (\ref{19}), we obtain
\begin{equation}\label{22}
\nabla ^{L,C} _{J\theta ^{\#}} =
-\frac{1}{2}(J\theta \circ \Box + \Box \circ J\theta ) + \frac{1}{2} dJ\theta -
\frac{1}{4}\iota _{J\theta ^{\#}} d^c \Omega .
\end{equation}
Hence, (\ref{15}) and (\ref{22}) yield
\begin{equation}\label{23}
\nabla ^{B,C} _{J\theta ^{\#}} =
-\frac{1}{2}(J\theta \circ \Box + \Box \circ J\theta ) + \frac{1}{2} dJ\theta .
\end{equation}
Now, by (\ref{18}) and (\ref{23}), we obtain
\begin{eqnarray}\label{24}
\| \nabla ^{B,B} \psi \| ^2 = & & \| \nabla ^{B,C} \psi \| ^2 -
(n-1)\Re (iJ\theta \Box  \psi ,\psi ) + \frac{n-1}{2}i(dJ\theta \psi ,\psi )
\\& & + \frac{(n-1)^2}{4}(|\theta |^2 \psi ,\psi). \nonumber
\end{eqnarray}
Hence, it follows from (\ref{17}) and (\ref{24}) that
\begin{eqnarray}\label{25}
\| \Box \psi \| ^2 = & & \| \nabla ^{B,B} \psi \| ^2 +
(n-1)\Re (iJ\theta \Box  \psi ,\psi ) + (\frac{s}{4} \psi ,\psi) +
\frac{i}{2}(\rho ^C \psi ,\psi ) \\& &+ \frac{1}{4}(dd^c \Omega \psi , \psi ) -
\frac{1}{8}(|d^c \Omega |^2 \psi ,\psi). \nonumber
\end{eqnarray}
It is clear that $|d^c \Omega |^2 = |d \Omega |^2$. Thus, using (\ref{25}),
(\ref{6}), (\ref{9}) and the fact that $\nabla ^{B,B}$ coincides with
$\nabla ^B$ restricted to $\Lambda ^{0,\bullet }M$, we obtain (\ref{16}).
\hfill {\bf Q.E.D.}

\vspace{3mm}
We recall that a $(1,1)$-form $\alpha $ is said to be positive (resp.
non-negative) if the corresponding symmetric tensor $A(X,Y)=\alpha (JX,Y)$ is
positive (resp. non-negative).

The following algebraic lemma is a direct consequence of the proof of
Theorem~1.1 and Remark~2.1.3 in \cite{H}.

\begin{lem}\label{lem2}
Let $\alpha $ be a $(1,1)$-form and let $a$ be its trace. If $\alpha $ is
positive (resp. non-negative), then as endomorphisms of
$\Lambda ^{0,p}$ $\frac{a}{2}Id + i\alpha $ is
positive definite (resp. non-negative definite) for $p>0$ and $-\frac{a}{2}Id + i\alpha $ is
negative definite (resp. non-positive definite) for $p<n$.
\end{lem}
Of course, $\frac{a}{2}Id + i\alpha $ is zero on $\Lambda ^{0,0}$ and $-\frac{a}{2}Id + i\alpha $ is zero on $\Lambda ^{0,n}$.

\subsection{Proof of Theorem~\ref{th2}}

Since $dd^c \Omega =0$, it follows
from (\ref{10}) that $|d\Omega |^2 = (n-1)^2 |\theta |^2 + (n-1)d^* \theta $.
The condition $(dd^c)^* \Omega =0$ is equivalent to $d^* \theta =0$.
The $(2,0)$- and (0,2)-forms send $\Lambda ^{0,p}M$
into $\Lambda ^{0,p-2}M$ and $\Lambda ^{0,p+2}M$ respectively.
Thus, if $\psi _p \in \Gamma (\Lambda ^{0,p}M)$, then
$(\rho ^B \psi _p, \psi _p) = ((\rho ^B)^{(1,1)} \psi _p, \psi _p)$, where
$(\rho ^B)^{(1,1)}$ is the $(1,1)$-part of $\rho ^B$. Hence, if
$\Box \psi _p =0$, Lemma~\ref{lem1} yields
$$0 = \| \nabla ^B \psi _p\| ^2 + (\frac{b}{4} \psi _p,\psi _p) +
\frac{i}{2}((\rho ^B)^{(1,1)} \psi _p,\psi _p).$$
The trace of a 2-form is equal to the trace of its $(1,1)$-part. Thus, the
assertions of the theorem follow from the latter equality and
Lemma~\ref{lem2}  \hfill {\bf Q.E.D.}

\vspace{3mm}
Let $h^{q,p}=\dim H^p (M,{\mathcal O}(\Lambda^{q,0} M))$ be the Hodge numbers of $(M,J)$. In particular, $h^{0,p}=\dim H^p (M,{\mathcal O})$

\begin{co}
Let $(M,g,J)$ be as in Theorem~\ref{th2}. Then the dimension of
the isometry group is greater or equal to $2h^{0,1}$.
\end{co}
\noindent {\it Proof:} Formula (\ref{3}) shows that every 1-form parallel with
respect to
the Bismut connection generates a Killing vector field. Thus the result follows
from Theorem~\ref{th2}. \hfill {\bf Q.E.D.}

\vspace{3mm}
The following lemma will be used in the proofs of Theorem~\ref{th3} and Theorem~\ref{th4}.

\begin{lem}\label{pro1}
Let $(M,g,J)$ be a $2n$-dimensional ($n>1$) Hermitian manifold such that the canonical Weyl connection preserves the complex structure. Let
$r ^W (X,Y) = Ric ^W (X,JY)$. Then $r ^W$ is a $(1,1)$-form. Further,

a) If $(M,g,J)$ is locally conformally K\"ahler, then
\begin{equation}\label{102}
\rho ^B = r^W - \frac{n-2}{2}d(J\theta).
\end{equation}

b) If $(M,g,J)$ is $4$-dimensional, then
\begin{equation}\label{103}
\rho ^B = r^W - (dJ\theta)^{(2,0)+(0,2)},
\end{equation}
where $(dJ\theta)^{(2,0)+(0,2)}$ denotes the $(2,0)+(0,2)$-part of $dJ\theta$.
\end{lem}

\noindent {\it Proof:} If $(M,g,J)$ is locally conformally K\"ahler, then $\nabla ^W$ locally is the Levi-Civita connection of a K\"ahler metric and hence $r ^W$ is a $(1,1)$-form. If $(M,g,J)$ is $4$-dimensional, the same is shown in \cite{GI}, formula (14).
\\
Since $\nabla ^W$ preserves the complex structure, it gives rise to a connection on $K^{-1}$. The equality (\ref{101}) is equivalent to
$d^c \Omega = J\theta \land \Omega $, so it follows from (\ref{3}) and (\ref{4}) that
\begin{equation}\label{104}
\nabla ^B _X \varphi = \nabla ^W _X \varphi + ( - \frac{n-2}{2}iJ\theta (X) + \frac{n}{2}\theta (X)) \varphi, \qquad
\varphi \in \Gamma (K^{-1}).
\end{equation}
Let $i\rho ^W$ be the curvature 2-form of $\nabla ^W$ as a connection on $K^{-1}$. Then (\ref{104}) implies
\begin{equation}\label{105}
i\rho ^B = i\rho ^W - \frac{n-2}{2}idJ\theta + \frac{n}{2}d\theta .
\end{equation}
Now, if $(M,g,J)$ is locally conformally K\"ahler, then $d\theta =0$. Since $\nabla ^W$ is locally the Levi-Civita connection of a K\"ahler metric, $\rho^W = r^W$. Thus (\ref{102}) follows from (\ref{105}).
\\
If $(M,g,J)$ is $4$-dimensional, then a direct computation shows
that
$$i\rho ^W = ir^W - d\theta - i(dJ\theta)^{(2,0)+(0,2)}.$$
By this equality and (\ref{105}) we obtain (\ref{103}).
\hfill {\bf Q.E.D.}

\subsection{Proof of Theorem~\ref{th4}}

Since $Ric^W$ is conformally
invariant, we can choose in the conformal class $c$ the Gauduchon metric
$g$ with respect
to which $d^* \theta =0$, or equivalently $(dd^c)^* \Omega =0$. It follows from (\ref{103}) that $(\rho^B)^{(1,1)} = r^W$ and therefore $(\rho^B)^{(1,1)}$ is non-negative everywhere on $M$ and strictly positive at
some
point. On 4-dimensional manifolds the condition $(dd^c)^* \Omega =0$ is equivalent to  $dd^c \Omega =0$. Hence, we can apply Theorem~\ref{th2} to obtain that $H^p (M,{\mathcal O}) =0$, $p=1,2$. The fact that $(\rho^B)^{(1,1)} = r^W$ implies also that $b \ge 0$ and $b$ is not identically zero. Thus, by Proposition~\ref{tpr1} all the plurigenera
of $(M,J)$ vanish. So, we have that $h^{0,1} = 0$ and $p_2 = 0$ and by the
Castelnuovo  criterion (cf. \cite{BPV}) $(M,J)$ must be a rational surface. The positivity of $c_1 ^2$ also follows from the fact that $(\rho^B)^{(1,1)}$ is non-negative everywhere and strictly positive at
some
point. \hfill {\bf Q.E.D.}

\vspace{3mm}
\noindent {\bf Remark~2}. From Lemma~\ref{pro1} and the main result in
\cite{GI}
we deduce that if on a compact Hermitian surface $(M,g,J)$ the $(1,1)$-part of
the
Ricci form of the Bismut connection is a scalar multiple of the K\"ahler form
at every point on $M$, then $(M,g,J)$ is conformally equivalent either to a
K\"ahler Einstein surface or to a Hopf surface.

\subsection{Proof of Theorem~\ref{th3}}

Let $\alpha $ be a 2-form. We have
$\alpha \land \Omega = \alpha \circ \Omega+ <\alpha , \Omega > + \varphi $
as endomorphisms of $\Lambda ^{0,\bullet }M$, where $\varphi $ is
$(2,0)+(0,2)$-form, $\varphi = 2i \alpha ^{(2,0)} - 2i \alpha ^{(0,2)}$. Hence,
if $\psi _p \in \Gamma (\Lambda ^{0,p}M)$, then
$$<\alpha \land \Omega \psi _p , \psi _p>= <\alpha \Omega \psi_p, \psi _p> +
<\alpha , \Omega > <\psi_p, \psi _p>.$$
Since $\Omega $ acts on $\Lambda ^{0,p}M$ as multiplication by $(n-2p)i$, we
obtain
\begin{equation}\label{25.5}
<\alpha \land \Omega \psi _p , \psi _p>= (n-2p)i<\alpha \psi_p, \psi _p> +
<\alpha , \Omega > <\psi_p, \psi _p>.
\end{equation}
Now, if $d\Omega = \theta \land \Omega $, then
$dd^c \Omega = (dJ\theta - J\theta \land \theta ) \land \Omega $. Since
$<J\theta \land \theta ,\Omega >= |\theta |^2$, using (\ref{8.5}) and
(\ref{25.5}) we get
$$<dd^c \Omega \psi_p, \psi _p>=
(n-2p)i<(dJ\theta - J\theta \land \theta ) \psi_p, \psi _p> +
((n-2)|\theta |^2 + d^* \theta ) <\psi_p, \psi _p>.$$
Thus, if $d^* \theta =0$ and $\Box \psi _p =0$, Lemma~\ref{lem1} and (\ref{11})
yield
\begin{eqnarray}\label{426}
0 = & & \| \nabla ^B \psi _p\| ^2 + \frac{1}{4}
((b+n(n-2)|\theta |^2)\psi _p,\psi _p) +
\frac{i}{2}(\rho ^B \psi _p,\psi _p) \\
& & +
\frac{n-2p}{4}i((dJ\theta - J\theta \land \theta ) \psi_p, \psi _p). \nonumber
\end{eqnarray}
Hence, when $n>2$ is even and $p=\frac{n}{2}$,   by
Lemma~\ref{lem2} we obtain  a).
\\
Suppose now that $(M,g,J)$ is a generalized Hopf manifold. In this case
\begin{equation}\label{427}
d(J\theta)=|\theta|^2 \Omega + \theta \wedge J\theta
\end{equation}
(see formula (2.8) in \cite{V4}). Thus
\begin{equation}\label{428}
i((dJ\theta - J\theta \land \theta ) \psi_p, \psi _p) = -(n-2p)(|\theta |^2 \psi_p, \psi _p) - 2i((J\theta \land \theta ) \psi_p, \psi _p).
\end{equation}
It follows from Lemma~\ref{lem2} that for $0<p<n$
$$-(|\theta |^2 \psi_p, \psi _p) \le i((J\theta \land \theta ) \psi_p, \psi _p) \le (|\theta |^2 \psi_p, \psi _p).$$
Hence, using (\ref{428}) we see that if $n-2p \ge 0$, then
\begin{equation}\label{429}
\frac{n-2p}{4}i((dJ\theta - J\theta \land \theta ) \psi_p, \psi _p) \ge  -\frac{(n-2p)(n-2p+2)}{4}(|\theta |^2 \psi_p, \psi _p)
\end{equation}
and if $n-2p \le 0$, then
\begin{equation}\label{430}
\frac{n-2p}{4}i((dJ\theta - J\theta \land \theta ) \psi_p, \psi _p) \ge  -\frac{(n-2p)(n-2p-2)}{4}(|\theta |^2 \psi_p, \psi _p) .
\end{equation}
Substituting (\ref{429}) and (\ref{430}) in (\ref{426}) we obtain
\begin{equation}\label{aa3}
0 \ge \| \nabla ^B \psi _p\| ^2 + \frac{1}{4}
(b+4(n-p)(p-1)|\theta |^2)\psi _p,\psi _p) +
\frac{i}{2}(\rho ^B \psi _p,\psi _p), \qquad n-2p \ge 0,
\end{equation}
$$0 \ge \| \nabla ^B \psi _p\| ^2 + \frac{1}{4}
(b+4p(n-p-1)|\theta |^2)\psi _p,\psi _p) +
\frac{i}{2}(\rho ^B \psi _p,\psi _p), \qquad n-2p \le 0,$$
Applying Lemma~\ref{lem2} to $(\rho^B)^{(1,1)}$ and using that $\theta \not =0$,
we get that $\psi _p =0$, i.e. $h^{0,p} =0$ for $1<p<n-1$. Applying Proposition~\ref{tpr1}, we get $p_1=0$. Hence $h^{0,n}=0$ and by Serre duality $h^{n,0}=0$.

Further, let $\psi$ be a $\bar{\partial}$-harmonic (0,1)-form. Then $\nabla^B
\psi=0$ by (\ref{aa3}). Using  $d\Omega=\theta \wedge
\Omega$, we obtain
$$
0=\bar{\partial}\psi = \theta ^{(0,1)}\wedge \psi,
$$
where $\theta ^{(0,1)}$ denotes the (0,1)-part of $\theta$. Hence $\psi
=f\theta^{(0,1)}$, where $f$ is a smooth function globally defined on $M$.
It is easy to see using (\ref{3}), $d\Omega =\theta \wedge \Omega$ and
$\nabla^L\theta =0$ that $\nabla^B\theta=0$ and that $\theta^{(0,1)}$ is
$\bar{\partial}$-harmonic. Thus, $\psi=f\theta^{(0,1)}$ leads to $f=const$ and
$h^{0,1}=1$. The vanishing of $h^{0,n-1}$ and $h^{p,0}$, $p=1,2,\ldots,n$ is an easy consequence
from the results of Tsukada \cite{Ts}, who proved that the Hodge numbers of every
generalized Hopf manifold satisfy the following relations:
\begin{equation}\label{hop}
h^{n,0}=h^{n-1,0}, \quad
h^{0,p}=h^{p,0}+h^{p-1,0}, \quad p\le n-1,
\end{equation}
\begin{equation}\label{hop1}
h^{1,0}=\frac{1}{2}(b_1(M)-1),
\quad h^{0,1}=\frac{1}{2}(b_1(M)+1).
\end{equation}
We also get $b_1=1$ from (\ref{hop1}).

Another proof of b) in Theorem~\ref{th3} can be obtained as follows. By (\ref{102}) and (\ref{426}) we have that
$$r^W = \rho^B + \frac{n-2}{2}(|\theta|^2 \Omega + \theta \wedge J\theta ).$$
Hence, $r^W$ is a non-negative $(1,1)$-form, i.e. $Ric^W$ is non-negative, and b) is proved applying Theorem~1.3 of \cite{AI}.  \hfill {\bf Q.E.D.}

\section{Examples}
\subsection{Compact Lie groups}\label{sec5.1}

Let $G$ be a compact connected Lie group with Lie algebra $\mathfrak{g}$. Denote by $\mathfrak{g}^c$ the complexification of $\mathfrak{g}$. Any left-invariant almost complex structure $J$ on $G$ is determined by its restriction on $\mathfrak{g}$, or equivalently, by the subspace $\mathfrak{s} \subset \mathfrak{g}^c$ of $(1,0)$-vectors. It is clear that
$$\mathfrak{s} \cap \mathfrak{g} = \{0\}, \qquad \mathfrak{g}^c = \mathfrak{s} \oplus \overline{\mathfrak{s}} $$
and $J$ is integrable iff $[\mathfrak{s},\mathfrak{s}] \subset \mathfrak{s}$, i.e. when $\mathfrak{s}$ is a complex Lie subalgebra of $\mathfrak{g}^c$. Such a subalgebra $\mathfrak{s}$ is called a Samelson subalgebra \cite{pit}. Samelson \cite{Sam} first constructed examples of left-invariant complex structures on compact Lie groups.

The construction is as follows. Let $T$ be a maximal torus in $G$, $\mathfrak{t}$ its Lie algebra and $\mathfrak{t}^c$ the complexification of $\mathfrak{t}$. Suppose a set $\alpha_1,\dots ,\alpha_m \in \mathfrak{t}^*$ of positive roots is chosen. Then all the roots are $\pm \alpha_1,\dots ,\pm \alpha_m$ and we have the $ad(T)$-invariant decomposition
\begin{equation}\label{551}
\mathfrak{g}^c = \mathfrak{t}^c \oplus \sum_{j=1}^m \mathfrak{s}_{\alpha_j} \oplus \sum_{j=1}^m \mathfrak{s}_{-\alpha_j},
\end{equation}
where
\begin{equation}\label{552}
\mathfrak{s}_{\pm \alpha_j} = \{Z \in \mathfrak{g}^c : [X,Z]= \pm 2\pi i\alpha_j(X)Z \quad  \forall X\in \mathfrak{t}\}.
\end{equation}
Now choose an almost complex structure on $\mathfrak{t}$, i.e. a subspace $\mathfrak{a} \subset \mathfrak{t}^c$ such that
$$\mathfrak{a} \cap \mathfrak{t} = \{0\}, \qquad \mathfrak{a} \oplus \overline{\mathfrak{a}} = \mathfrak{t}^c.$$
Then it is clear that
\begin{equation}\label{553}
\mathfrak{s} = \mathfrak{a} \oplus \sum_{j=1}^m \mathfrak{s}_{\alpha_j}
\end{equation}
is a Samelson subalgebra of $\mathfrak{g}^c$ and hence gives rise to a left-invariant complex structure on $G$.

It is proved by Pittie \cite{pit} that, conversely, any left invariant complex structure on a compact Lie group can be obtained as above. We sketch a proof of this fact. Let $\mathfrak{s}$ be the Samelson subalgebra corresponding to a left invariant complex structure $J$ on $G$. Define
$$\mathfrak{t}= \{X \in \mathfrak{g} : ad(X)\mathfrak{s} \subset \mathfrak{s} \},$$
i.e. $\mathfrak{t}$ consists of all elements of $\mathfrak{g}$ which preserve the decomposition $\mathfrak{g}^c = \mathfrak{s} \oplus \overline{\mathfrak{s}}$. It is easy to see that $\mathfrak{t}$ is a $J$-invariant subalgebra of $\mathfrak{g}$. Hence $\mathfrak{t}$ is a complex Lie algebra. Let $T$ be the closed connected subgroup of $G$ which corresponds to $\mathfrak{t}$. Then $T$ is a compact complex Lie group and hence a torus. In particular, $\mathfrak{t}$ is abelian. We can choose an $ad(T)$-invariant inner product on $\mathfrak{g}$, which is also $J$-invariant. Then with respect to this inner product we have an orthogonal $ad(T)$-invariant decomposition as in (\ref{553}), where $\alpha_j \in \mathfrak{t}^*$, $\mathfrak{s}_{\alpha_j}$ is defined as in (\ref{552}) and
$$\mathfrak{a} = \{Z \in \mathfrak{s} : [X,Z]=0 \quad \forall X\in \mathfrak{t}\}.$$
It follows easily from the definition of $\mathfrak{t}$ that $\mathfrak{t}^c = \mathfrak{a} \oplus \overline{\mathfrak{a}}$, i.e. $\mathfrak{t}$ is a maximal abelian subalgebra of $\mathfrak{g}$ and $T$ is a maximal torus. Thus (\ref{551}) is satisfied and hence $\pm \alpha_1,\dots ,\pm \alpha_m$ are all the roots. It follows from $[\mathfrak{s},\mathfrak{s}] \subset \mathfrak{s}$ that if $\alpha_j + \alpha_k$ is a root, then  $\alpha_j + \alpha_k = \alpha_l$ for some $l$. Thus we can take the set of positive roots to be $\{ \alpha_1,\dots ,\alpha_m \}$. Hence, any left-invariant complex structure on a compact Lie group is determined by a choice of a maximal torus, a complex structure on its Lie algebra and a choice of positive roots.

In \cite{pit} Pittie has described the moduli of left-invariant complex structures on compact Lie groups. For semi-simple groups he has also studied the Dolbeault cohomology rings of such complex structures using Bott's Lie algebraic description of these rings \cite{Bot}. Among other results Pittie has proved that the Hodge numbers $h^{0,p}$ of a left-invariant complex structure on a semi-simple group $G$ are the same as the corresponding Hodge numbers of the maximal torus $T$. It is not hard to see that this is true for any compact Lie group without requiring semi-simplicity: It is clear that $G$ is a holomorphic principal bundle over $G/T$ with structure group $T$. It is also well-known that $G/T$ has positive first Chern class (see e.g. \cite{Bes}, chapter 8) and hence by Bochner-Kodaira vanishing theorem $h^{0,p}(G/T) = 0$ for $p>0$. Thus, applying the results of Borel ({\bf 7.5} of \cite{bor}) to the bundle $G(G/T,T)$ we obtain that $h^{0,p}(G) \le h^{0,p}(T)$. On the other hand, it is clear that $h^{0,p}(G) \ge h^{0,p}(T)$ and hence $h^{0,p}(G) = h^{0,p}(T)$.

Now suppose that $G$ is a compact Lie group with a left-invariant complex structure $J$ which is compatible with a bi-invariant metric $g$. In this case the left-invariant connection on $G$ is a Hermitian connection on $(G,g,J)$ and because of the bi-invariance of the metric its torsion tensor is totally skew-symmetric. This means that the left-invariant connection is the Bismut connection of the given Hermitian structure. Since the left-invariant connection is flat, then $\rho^B =0$. It is easy to see also that $dd^c \Omega =0$ (i.e. the torsion 3-form is closed) and $(dd^c)^* \Omega =0$ (i.e. $d^* \theta =0$). Thus $(G,g,J)$ is a manifold which satisfies the assumptions of Theorem~\ref{th2}. In this case the above mentioned result for the Hodge numbers $h^{0,p}$ can be obtained as a corollary of that theorem.

\begin{co}\label{prop51}
Let $G$ be a compact Lie group with $rank \, G =2r$ and $J$ be a left-invariant complex structure on $G$ which is compatible with a bi-invariant metric $g$. Then the Hodge numbers of $(G,J)$ are $h^{0,p} = {r \choose p}$
\end{co}

\noindent {\it Proof:} We keep the above introduced notations. It follows from a) of Theorem~\ref{th2} that every $\overline{\partial}$-harmonic $(0,p)$-form on $(G,g,J)$ is left-invariant. Let $Z_1, \dots ,Z_{m+r}$ be an orthonormal basis of $\mathfrak{s}$ (i.e. $g(Z_j, \overline{Z}_k)= \delta_{jk}$) such that $Z_j \in \mathfrak{s}_{\alpha_j}$ for $j=1,\dots,m$ and $Z_j \in \mathfrak{a}$ for $j>m$. Denote by $\zeta^1,\dots,\zeta^{m+r}$ the dual basis. Let $\overline{p}: \mathfrak{g}^c \longrightarrow  \overline{\mathfrak{s}}$ be the projection. Then for a left-invariant $(0,p)$-form $\varphi$ and $\overline{W}_1,\dots,\overline{W}_{p+1} \in \overline{\mathfrak{s}}$ we have
\begin{equation}\label{504}
\overline{\partial}\varphi(\overline{W}_1,\dots,\overline{W}_{p+1}) = \sum_{j<k} (-1)^{j+k} \varphi([\overline{W}_j,\overline{W}_k],\overline{W}_1,\dots,\widehat{\overline{W}}_j,\dots,\widehat{\overline{W}}_k,\dots,\overline{W}_{p+1}),
\end{equation}
\begin{equation}\label{505}
\overline{\partial}^* \varphi(\overline{W}_1,\dots,\overline{W}_{p-1}) = \sum_{j=1}^{p-1} \sum_{k=1}^{m+r} (-1)^{j-1} \varphi(\overline{Z}_k,\overline{p}([Z_k,\overline{W}_j]),\overline{W}_1,\dots,\widehat{\overline{W}}_j,\dots,\overline{W}_{p-1}),
\end{equation}
where $\widehat{\overline{W}}_j$ means that $\overline{W}_j$ has to be deleted. Let $\mathfrak{n}= \sum_{j=1}^m \mathfrak{s}_{\alpha_j}$. Hence $\overline{\mathfrak{s}} = \overline{\mathfrak{n}} \oplus \overline{\mathfrak{a}}$ and
$$\Lambda^p \overline{\mathfrak{s}}^* = \sum_{d=0}^p \Lambda^d\overline{\mathfrak{n}}^* \otimes \Lambda^{p-d}\overline{\mathfrak{a}}^*.$$
All the forms in $\Lambda^p\overline{\mathfrak{a}}^*$ are $\overline{\partial}$-harmonic and $\dim_{\bf C} \Lambda^p\overline{\mathfrak{a}}^* = {r \choose p}$. Hence, it remains to show that there are no $\overline{\partial}$-harmonic forms in $\sum_{d=1}^p \Lambda^d\overline{\mathfrak{n}}^* \otimes \Lambda^{p-d}\overline{\mathfrak{a}}^*$. Using (\ref{504}) and (\ref{505}) this can be proved by induction. We assume that a $\overline{\partial}$-harmonic form $\varphi$ belongs to $\sum_{d=q}^p \Lambda^d\overline{\mathfrak{n}}^* \otimes \Lambda^{p-d}\overline{\mathfrak{a}}^*$ and the $\Lambda^q \overline{\mathfrak{n}}^* \otimes \Lambda^{p-q} \overline{\mathfrak{a}}^*$-part of $\varphi$ is of the form
$$\sum_{i_1 < \dots <i_q < m} \overline{\zeta}^{i_1} \wedge \dots \wedge \overline{\zeta}^{i_q} \wedge \overline{p}(\alpha_{i_1} + \dots + \alpha_{i_q}) \wedge \psi_{i_1 \dots i_q},$$
where $\psi_{i_1 \dots i_q} \in \Lambda^{p-q-1}\overline{\mathfrak{a}}^*$. Then we prove that the $\Lambda^q \overline{\mathfrak{n}}^* \otimes \Lambda^{p-q} \overline{\mathfrak{a}}^*$-part of $\varphi$ is zero and the $\Lambda^{q+1} \overline{\mathfrak{n}}^* \otimes \Lambda^{p-q-1} \overline{\mathfrak{a}}^*$-part of $\varphi$ is of the form
$$\sum_{i_1 < \dots <i_{q+1} < m} \overline{\zeta}^{i_1} \wedge \dots \wedge \overline{\zeta}^{i_{q+1}} \wedge \overline{p}(\alpha_{i_1} + \dots + \alpha_{i_{q+1}}) \wedge \psi_{i_1 \dots i_{q+1}},$$
where $\psi_{i_1 \dots i_{q+1}} \in \Lambda^{p-q-2}\overline{\mathfrak{a}}^*$. \hfill
{\bf Q.E.D.}

\subsection{Hermitian manifolds with {\boldmath $(\rho ^B)^{(1,1)} >0$}}\label{sec5.2}

Let $(M,\widetilde{g},J)$ be a compact 4-dimensional K\"ahler manifold with positive Ricci form. Let $\alpha$ be a $(1,1)$-form and let $\Omega = \widetilde{\Omega} + \alpha$ ($\widetilde{\Omega}$ is the K\"ahler form of $(M,\widetilde{g},J)$). Clearly, if $\alpha$ has sufficiently small $C^2$-norm, then $\Omega$ will be the K\"ahler form of a Hermitian metric $g$ and the $(1,1)$-part of the Ricci form of the Bismut connection of $(M,g,J)$ will be positive. Let $g_0 \in [g]$ be the Gauduchon metric. Because of the conformal invariance of $(\rho ^B)^{(1,1)}$, the $(1,1)$-part of the Ricci form of the Bismut connection of $(M,g_0,J)$ will also be positive. For generic $\alpha$ the manifold  $(M,g_0,J)$ is not K\"ahler. So, we obtain non-K\"ahler 4-dimensional manifolds which satisfy the assumptions of Theorem~\ref{th2} b) (and hence also manifolds satisfying the assumptions of Theorem~\ref{th4} with a conformal class of Hermitian metrics not containing a K\"ahler metric).

Now let $(M',g',J')$ be a compact $(2n-4)$-dimensional K\"ahler manifold with positive Ricci form and $(M'',g'',J'')$ be a compact 4-dimensional Hermitian non-K\"ahler manifold with positive $(1,1)$-part of the Ricci form of the Bismut connection. Then it is easy to verify that $M= M' \times M''$ with the product Hermitian structure is a compact $2n$-dimensional Hermitian non-K\"ahler manifold with $dd^c$-harmonic K\"ahler form and positive $(1,1)$-part of the Ricci form of the Bismut connection. Thus there are non-K\"ahler manifolds satisfying the assumptions of Theorem~\ref{th2} b) in any dimension.

\bigskip {\bf Authors' address:}\\[2mm] Bogdan Alexandrov,
Stefan Ivanov\\ University of Sofia, Faculty of
Mathematics and Informatics, Department of Geometry, \\ 5 James
Bourchier Blvd, 1126 Sofia, BULGARIA.\\ E-mail: B.Alexandrov:
{\tt alexandrovbt@fmi.uni-sofia.bg} \par
S.Ivanov: {\tt ivanovsp@fmi.uni-sofia.bg}
\end{document}